\documentclass[a4paper,12pt]{article}
\usepackage{amssymb}
\usepackage{latexsym}
\usepackage{amsmath}
\usepackage{amsthm}
\usepackage{amsfonts}
\usepackage{amssymb}

\usepackage{pstricks}

\newcommand{\rec}[1]{{(\ref{#1})}} 
\newtheorem{Th}{Theorem}[section]
\newtheorem{Prop}{Proposition}[section]

\def\Xint#1{\mathchoice
{\XXint\displaystyle\textstyle{#1}}
{\XXint\textstyle\scriptstyle{#1}}
{\XXint\scriptstyle\scriptscriptstyle{#1}}
{\XXint\scriptscriptstyle\scriptscriptstyle{#1}}
\!\int}
\def\XXint#1#2#3{{\setbox0=\hbox{$#1{#2#3}{\int}$ }
\vcenter{\hbox{$#2#3$ }}\kern-.6\wd0}}
\newcommand{\R}{\mathbb{R}}
\newcommand{\C}{\mathbb{C}}

\newcommand{\dd}{\mathrm{d}} 

\DeclareMathOperator{\divop}{div}

\def\dashint{\Xint-}
\def\Om{\Omega}

\def\p{\partial}

\newcommand{\be}{\begin{equation}}
\newcommand{\ee}{\end{equation}}
\newcommand{\bes}{\begin{equation*}}
\newcommand{\ees}{\end{equation*}}

\def\ti{\tilde}
\def\lf{\left}
\def\rg{\right}

\def\la{\lambda}

\def\ds{\displaystyle}

\def\Om{\Omega}

\def\p{\partial}

\definecolor{babyblue}{rgb}{0.54, 0.81, 0.94}
\definecolor{bittersweet}{rgb}{1.0, 0.44, 0.37}
\definecolor{brightmaroon}{rgb}{0.76, 0.13, 0.28}
\definecolor{icterine}{rgb}{0.99, 0.97, 0.37}
\definecolor{indiagreen}{rgb}{0.07, 0.53, 0.03}
\definecolor{aquamarine}{rgb}{0.5, 1.0, 0.83}

\definecolor{amaranth}{rgb}{0.9, 0.17,0.31}
\definecolor{amethyst}{rgb}{0.6, 0.4, 0.8}
\definecolor{bostonuniversityred}{rgb}{0.8, 0.0, 0.0}
\definecolor{babyblue}{rgb}{0.54, 0.81, 0.94}
\definecolor{coolblack}{rgb}{0.0, 0.18, 0.39}
\definecolor{cobalt}{rgb}{0.0, 0.28, 0.67}
\definecolor{amber}{rgb}{1.0, 0.75, 0.0}
\definecolor{azure(colorwheel)}{rgb}{0.0, 0.5, 1.0}
\definecolor{shamrockgreen}{rgb}{0.0, 0.62, 0.38}
\definecolor{brinkpink}{rgb}{0.98, 0.38, 0.5}

\begin{document}
\title{\bf  Remarks  on Neumann boundary problems involving Jacobians   }
\author{   Francesca Da Lio\thanks{Department of Mathematics, ETH Z\"urich, R\"amistrasse 101, 8092 Z\"urich, Switzerland.}  
\and Francesco Palmurella$^*$  }
\maketitle 
\begin{abstract}In this short note we explore the validity of Wente-type estimates for Neumann boundary problems involving Jacobians.
We  show in particular that such estimates do not in general hold under the same hypotheses  on the data for Dirichlet boundary problems. \end{abstract}
\par
{\small { \bf  MSC 2010.  35J25, 30C70, 35B65}
\section{Introduction}
Integrability by compensation has played a central role in the last decades  in the geometric analysis of conformally invariant problems. At the center
of this theory there is the Wente's discovery \cite{Wente} that the distribution:
$$
\varphi(x)=\log|x|\star [\partial_{x_1}a\partial_{x_2}b-\partial_{x_2}a\partial_{x_1}b]
$$
with $\nabla a,\nabla b\in L^2(\R^2)$ is in $(L^{\infty}\cap W^{1,2})(\R^2)$ and the following estimate holds true:
$$
\|\varphi\|_{L^{\infty}(\R^2)}\le C\|\nabla a\|_{L^2(\R^2)}\|\nabla b\|_{L^2(\R^2)}.$$
It has been observed by Brezis and Coron in \cite{BrCo} that a similar estimate holds also  if we consider the following Dirichlet problem:
\begin{Th}(\cite{Wente})\label{WenteTh}
Let $\Omega=\{(x,y)\in\R^2:x^2+y^2<1\}$ and let
$a,b\in H^1(\Omega)$. Then   the solution $u\in W^{1,1}_0(\Omega)$ to the problem:
\begin{equation}
\label{We}
\lf\{
\begin{array}{cc}
\ds-\Delta u=\p_{x_1}a\,\p_{x_2} b-\p_{x_2}a\,\p_{x_1} b&\quad\quad\mbox{ in }\Omega,\\[5mm]
u=0&\quad\quad\mbox{ on }\p \Omega,
\end{array}
\rg.
\end{equation}
is a continuous function in $\overline \Omega$ and its gradient belongs to $L^2(\Omega).$ Moreover there exists a constant $C_0=C_0(\Omega)$ such that
\begin{equation}\label{West}
\|u\|_{L^{\infty}(\Omega)}+\|\nabla u\|_{L^2(\Omega)}\le C_0\|\nabla a\|_{L^2(\Omega)}\|\nabla a\|_{L^2(\Omega)}.\end{equation}
\end{Th}
Extensive investigation on the problem \eqref{We} and various generalisations has been 
conducted, remarkably in \cite{CLMS}.

The goal of the present work is to explore to which extent an inequality like \rec{West} holds or not if we replace the  Dirichlet boundary condition with a Neumann boundary condition.
\par Our first main result gives a negative answer for general $a$ and $b$.  We consider for simplicity the unit disk $D^2:=\{(x,y)\in\R^2:x^2+y^2<1\}$ and we denote by $\nu$  the unit outward normal vector to
$\partial D^2$.   Then:
\begin{Th}\label{NegAnsw1}
There are  $a,b\in (L^{\infty}\cap {H}^{1})(D^2)$  such that 
the solution $\varphi$ with zero mean value of: 
\begin{equation}\label{NPHIntr}\left\{\begin{array}{cc}
-\Delta \varphi=\nabla a\cdot\nabla^{\perp}b& ~~~\mbox{in $D^2$},\\
\partial_{\nu}\varphi=0& ~~~\mbox{on $\partial D^2$},
\end{array}\right.
\end{equation} 
is not  in $H^{1}(D^2).$
\end{Th}
\par  A Wente-type estimate holds for \rec{NPHIntr} if ${\rm Trace}(a)=0$ or ${\rm Trace}(b)=0$. We call this the case
{\em vanishing Jacobian at the boundary}. Precisely,   the following result holds.\par
\begin{Th}
\label{th-vanish}
Let $\Omega$ be a smooth bounded domain of ${\R}^2$. Let $a\in W^{1,2}_0(\Omega )$ and $b\in W^{1,2}(\Omega )$ and let $u$ be a solution of \rec{NPHIntr}.
Then $\nabla^2 u\in L^{1}(\Om)$ and one has:
\begin{equation}
\label{II.2}
\|\nabla^2 u\|_{L^{1}(\Om)}\le C(\Om)\ \|\nabla a\|_{L^2(\Om)}\ \|\nabla b\|_{L^2(\Om)}. 
\end{equation}
\end{Th}

From estimate \eqref{II.2}, by means of improved Sobolev embeddings (see e.g. \cite{Hel})
we deduce that the estimate \eqref{We} holds. 
Theorem \ref{th-vanish} has been used by Rivi\`ere in \cite{Riv} and the proof will be given in \cite{DLPR}.

In some applications such as for instance in the analysis of the  {\em Poisson problem for elastic plates} (\cite{DLPR})  the  following Neumann boundary problem  
appears in a natural way:
{  \begin{equation}\label{NPCompIntr}\left\{\begin{array}{cc}
-\Delta w=\nabla a\cdot\nabla^{\perp}b ~~~\mbox{in $D^2$,}\\
\partial_{\nu}w=-\ a\partial_{\tau} b ~~~\mbox{on $\partial D^2$.}
\end{array}\right.
\end{equation} }
We observe that $H^{1}(D^2)$-solutions of the problem \rec{NPCompIntr} are critical points of the following Lagrangian:
\begin{equation}\label{lagra1}
{\mathcal{L}}(u;a,b)=\frac{1}{2}\int_{D^2}|\nabla u+a(\nabla^{\perp}b) |^2 dy_1 dy_2
\end{equation}

We will refer to the problem \rec{NPCompIntr} as the case of {\em compatible Neumann boundary conditions}.
Also in the case of  \rec{NPCompIntr} the assumption  $a,b\in (L^{\infty}\cap {H}^1)(D^2)$ is not enough to guarantee the boundedness
of the solution in $\bar D^2.$ 
\begin{Th}\label{NegAnsw2}
There are  $a,b\in L^{\infty}(D^2)\cap {H}^{1}(D^2)$  such that 
the solution $\varphi$ with zero mean value of: 
\begin{equation}\label{NPHIntr}\left\{\begin{array}{cc}
-\Delta \varphi=\nabla a\cdot\nabla^{\perp}b& ~~~\mbox{in $D^2$},\\
\partial_{\nu}\varphi =-\ a\partial_{\tau} b& ~~~\mbox{on $\partial D^2$,}
\end{array}\right.
\end{equation} 
is not  in $L^{\infty}(\bar D^2).$
\end{Th}

The boundedness of the solution is however obtained if  we assume a bit more on the data $a,b$. More precisely we get the following result.
\begin{Th}[$L^{2,1}$-case]\label{main}
Let 
 $ \nabla b\in L^{(2,1)}(D^2)$,\footnote{We denote by  
$L^{2,1}(\R^n)$  the  space of measurable functions satisfying
$$\int_{0}^{+\infty}|\{x\in\R^n~: |f(x)|\ge\lambda\}|^{1/2} d\lambda<+\infty \,\,.$$}   $a\in L^{\infty}(D^2)\cap {H}^{1}(D^2)$    
and let    $w\in W^{1,1}(D^2)$ be the  solution with zero mean value
to \rec{NPCompIntr} . Then $\nabla w\in  L^{(2,1)}(D^2)$ with:
\begin{eqnarray} \|\nabla w\|_{L^{(2,1)}(D^2)}&\le &C\|  a\|_{L^{\infty}}\|\nabla b\|_{L^{(2,1)}}.\end{eqnarray}
In particular:
\begin{equation}  \|w \|_{L^{\infty}(\bar D^2)}\le C\|  a\|_{L^{\infty}}\|\nabla b\|_{L^{(2,1)}},\end{equation}
 \end{Th}
We observe that the assumption $ \nabla b\in L^{(2,1)}(D^2)$ is in particular satisfied if $b\in {W}^{2,1}(D^2)$, see e.g \cite{Hel}.
We remark that if we assume that also $ \nabla a\in L^{(2,1)}(D^2)$ and $\bar a=\dashint_{D^2} a(y) dy=0$  then $a\in L^{\infty}(D^2)$ and
$\|  a\|_{L^{\infty}}\le C \|\nabla a\|_{L^{(2,1)}}$. In this case we can estimate $\nabla w$ as follows:
\begin{eqnarray}
\|\nabla w\|_{L^{(2,1)}(D^2)}&\le &C\|\nabla a\|_{L^{(2,1)}}\|\nabla b\|_{L^{(2,1)}}.\end{eqnarray}

\par 
The paper is organized as follows. In Section 2 we prove Theorem \ref{main} and in Section 3 we prove Theorems \ref{NegAnsw1} and
\ref{NegAnsw2}.

\section{Proof of Theorem \ref{main}.}
{\bf Step 1.}
We start by observing that
we can formulate problem \rec{NPCompIntr} as follows:
{ \begin{equation}\label{pb2}\left\{\begin{array}{cc}
{\rm div}[\nabla w+ a \nabla^\perp b]=0 &~~~\mbox{in $D^2$},\\[5mm]
\partial_{\nu} w=- a\partial_{\tau}b  &~~~\mbox{in $\partial D^2$}.
\end{array}\right.
\end{equation} }
Therefore there exists $C\in W^{1,2}_{0}(D^2)$ such that:
$$
\nabla^{\perp} C=\nabla w+  a\,\nabla^\perp b .$$
Therefore $C$ solves:
{ \begin{equation}\label{pb3}\left\{\begin{array}{cc}
-\Delta C=-{\rm div}(a\,\nabla b ) ~~~\mbox{in $D^2$},\\[5mm]
\partial_{\tau}C=0 ~~~\mbox{in $\partial D^2$}.
\end{array}\right.
\end{equation} }
Since $C$ is determined up to a constant, we can reduce to study the following Dirichlet problem:
{ \begin{equation}\label{pb3}\left\{\begin{array}{cc}
-\Delta C=-{\rm div}(a\,\nabla b ) ~~~\mbox{in $D^2$},\\[5mm]
C=0 ~~~\mbox{in $\partial D^2$}.
\end{array}\right.
\end{equation} }

{\bf Step 2.}
In this step and in the following we use basic facts about 
the theory of Calder\'on-Zygmund operators and interpolation theory,
for which we refer to
\cite{Hel,Ste}.
We first assume $b\in W^{1,p}(D^2). $  Let us set $f= -a\,\nabla b  \in L^{p}(D^2)$, we have:
\begin{eqnarray*}
\|f\|_{L^{p}(D^2)}&\le& C \|\nabla b\|_{L^p(D^2)}\|a\|_{L^{\infty}(D^2)}.\\
\end{eqnarray*}
We denote by $\tilde f=f\chi_{D^2}$
its extension by 0 to $\R^2$. 
We write $C=C_1+C_2$ where:
\begin{align*}
C_1(x)=\left(
-\frac{1}{2\pi}\log|\cdot|* \divop\tilde{f}
\right)(x),\quad x\in\R^2,
\end{align*}
and $C_2=C-C_1$ which is the solution to:
{ \begin{equation}\label{pb3}\left\{\begin{array}{cc}
-\Delta C_2=0~~~\mbox{in $D^2$},\\[5mm]
C_2=-C_1~~~\mbox{in $\partial D^2$}.
\end{array}\right.
\end{equation} }
We have: 
\begin{eqnarray*}
\nabla C_1(x_1,x_2) 
&=&\frac{1}{2\pi}\int_{\R^2 }\tilde f(y)\left[\frac{y-x}{|y-x|^3}\right] \dd y.
\end{eqnarray*}

The function
$${\mathcal{K}}(x,y)=\frac{y-x}{|y-x|^3} $$
is a C-Z operator. 
Since  $\tilde f\in L^p(\R )$ for every $p>1$ we have 
$$T[\tilde f](x):= \frac{1}{2\pi}\int_{\R^2}\tilde f(y)\left[\frac{y-x}{|y-x|^3}\right] dy\in L^p(\R^2)$$
and   $$
\|T[\tilde f]\|_{L^p}\le C_p\|\tilde f\|_{L^p}.$$
(see e.g. \cite{Hel,Ste}).

As far as $C_2$ is concerned,
since $C_1\in W^{1-1/p,p} (\partial D^2),$ then $C_2\in W^{1,p}(D^2)$ and:
$$\|\nabla C_2\|_{L^p(D^2)}\le C_p\|C_1\|_{ W^{1-1/p,p}}\le C_p \|f\|_{L^p(D^2)}.$$
In particular we get:

$$\|\nabla C\|_{L^p(D^2)}\le   C_p \|f\|_{L^p(D^2)},$$
and therefore:
$$\|\nabla w\|_{L^p(D^2)}\le   C_p \|f\|_{L^p(D^2)}\le C \|\nabla b\|_{L^p(D^2)}\|a\|_{L^{\infty}(D^2)}.$$
We remind that if  $p$ belongs to a compact interval $I\subset (0,\infty)$, the constant $C_p$ is uniformly bounded. \par
Now we define
$${\frak{G}}_p(D^2):=\{X\in L^p(D^2,\R^m):~~\mbox{ curl$( X)$}\}=0.$$
\footnote{We recall that for $X\in L^p(D^2,\R^m)$, $ \mbox{ curl$( X)=(-X^i_{x_2}+X^{i}_{x_1})_i$, $1\le i\le m$.}$ If ${\rm curl} X=0$ then there exists a $b\in W^{1,p}(D^2)$ such that $\nabla b=X$.}
Therefore if we fix $a\in L^{\infty}(D^2)$, the operator
$\tilde T\colon {\frak{G}}_p(D^2) \to L^p(D^2),$ $\nabla b\mapsto\nabla w$ is continuous for each $p>1$.

{\bf Step 3.}
If $a\in L^{\infty}(D^2)$ and $ \nabla b\in  L^{(2,1)}(D^2)$ then $f\in L^{(2,1)}(D^2)$ with
\begin{eqnarray*}
\|f\|_{L^{(2,1)}(D^2)}&\le& C 
\|\nabla b\|_{L^{(2,1)}(D^2)} \|a\|_{L^{\infty}(D^2)}.
\end{eqnarray*}

By interpolation and the previous step,
we get that  $\nabla w\in L^{(2,1)}(D^2)$  with:
\begin{equation}\label{nablaw}
\|\nabla w\|_{L^{(2,1)(D^2)}}\le C\|f\|_{L^{(2,1)}(D^2)}\le C \| a\|_{L^{\infty}(D^2)}\|\nabla b\|_{L^{(2,1)}(D^2)}.\end{equation}

for some $C>0$. \par
If we suppose that $a\in L^{\infty}$,  $  \nabla b\in  L^{(2,1)}(D^2)$ then from \rec{nablaw} it follows that
\begin{equation}\label{nablaw2}
\|\nabla w\|_{L^{(2,1)}(D^2)}\le C \| a\|_{L^{\infty}(D^2)}\|\nabla b\|_{L^{(2,1)}(D^2)},\end{equation}
and we conclude.~~\hfill$\Box$
\section{ Proof of Theorems \ref{NegAnsw1} and \ref{NegAnsw2} }
In this Section we provides counter-examples to Wente-type estimates for solutions to \rec{NPCompIntr} and \rec{NPHIntr} even in the case $a,b\in (H^1\cap L^{\infty})(D^2)$.
\subsection{A representation formula with estimates}
\label{subsec:repform}
Because of the conformal invariance of the problem \rec{NPCompIntr} we can reduce to consider the  analogous problem in $\R^2_+:$
\begin{equation}\label{pbhs}\left\{\begin{array}{cc}
-\Delta w=\nabla a\cdot\nabla^{\perp}b ~~~\mbox{in $\R^2_+$}\\
\partial_{\nu}w=-\ a\partial_{\tau} b ~~~\mbox{in $\partial \R^2_+.$}
\end{array}\right.
\end{equation}  
The Green function associated to the Neumann problem in the half-plane
$\mathcal{G}:\R^2_+\times\R^2_+\to\R$  is the solution, for every $x\in\R^2_+$
of the problem:
\begin{align*}
\left\{
\begin{aligned}
-\Delta_y\mathcal{G}(x,\cdot) &=\delta_x  &&\text{in }\R^2_+,\\
\partial_{\nu_y}\mathcal{G}(x,\cdot)&=0 &&\text{in }\partial \R^2_+,
\end{aligned}
\right.
\end{align*}
given by:
$${\mathcal{G}}(x,y)=-\frac{1}{2\pi}\left\{\log(|x-y|)+\log(|y-\tilde x|)\right\},$$
where $x=(x_1,x_2)$, $y=(y_1,y_2)$, $\tilde x=(x_1,-x_2)$. 

We are going to consider
the solution $w$ to \eqref{pbhs} obtained through the representation formula:
\begin{eqnarray}\label{reprform}
w(x)&=&\int\int_{R^2_+}{\mathcal{G}}(x,y)(-\Delta w)  dy+\int_{\partial\R^2_+}{\mathcal{G}}(x,y)\partial_{\nu} w d\sigma(y)\nonumber\\
&=&-\frac{1}{2\pi}\int\int_{R^2_+ }\left\{\log(|x-y|)+\log(|y-\tilde x|)\right\}\nabla a\cdot\nabla^{\perp}bdy\\
&&-\frac{1}{\pi}\int_{-\infty}^{+\infty}\log((y_1-x_1)^2+x^2_2)^{1/2} a\partial_{y_1} bdy_1\nonumber
\end{eqnarray}
and deduce a representation formula for \emph{its trace} at the boundary
$\partial\R^2_+$.

{\bf Step 1:} We assume for the moment that $a,b$ are in $C^\infty_c(\R^2).$
We integrate by parts \rec{reprform} and get:
\begin{eqnarray}\label{reprform2}
w(x)&=&-\frac{1}{2\pi}\int\int_{\R^2_+ }{\rm div}\left(\left\{\log(|x-y|)+\log(|y-\tilde x|)\right\} a\cdot\nabla^{\perp}b\right)dy\nonumber\\
&+&\frac{1}{2\pi}\int\int_{\R^2_+ }\nabla(\left\{\log(|x-y|)+\log(|y-\tilde x|)\right\})(a\cdot\nabla^{\perp}b)dy\nonumber\\
&-&\frac{1}{\pi}\int_{-\infty}^{+\infty}\log((y_1-x_1)^2+x^2_2)^{1/2} [a\partial_{y_1} b ]dy_1\nonumber\\
&=&\frac{1}{2\pi}\int_{\partial \R^2_+ }\left\{\log(|x-y|)+\log(|y-\tilde x|)\right\} a\partial_{y_1}b \,d\sigma(y)\\
&+&\frac{1}{2\pi}\int\int_{\R^2_+ }\nabla(\left\{\log(|x-y|)+\log(|y-\tilde x|)\right\}) a\nabla^{\perp}b)dy\nonumber\\
&-&\frac{1}{\pi}\int_{-\infty}^{+\infty}\log((y_1-x_1)^2+x^2_2)^{1/2} a\partial_{y_1} b]dy_1\nonumber\\
&=& { \frac{1}{2\pi}\int\int_{\R^2_+ }\nabla(\left\{\log(|x-y|)+\log(|y-\tilde x|)\right\})a \nabla^{\perp}b dy}. \nonumber
\end{eqnarray}

If $x\in \partial \R^2_+$, then: 
\begin{eqnarray}\label{reprform4}
w(x_1,0)&=& \frac{1}{\pi}\int\int_{\R^2_+ }\nabla(\left\{\log((x_1-y_1)^2+y_2^2)^{1/2} \right\})\cdot  [a \nabla^{\perp}b ]dy.
\end{eqnarray}
By translation invariance we evaluate \rec{reprform4} at $(0,0)$ and use polar
coordinates.
For every $r>0$ we set
$$ a_r=\frac{a(r,\pi)+a(r,0)}{2}=\frac{a(x_1,0)+a(-x_1,0)}{2}.$$

We have
\begin{eqnarray}\label{bcterm}
w(0,0)&=&\frac{1}{\pi}\int\int_{R^2_+ }\nabla( \log|y| )\cdot  [a \nabla^{\perp}b ]dy\nonumber\\
&=&-\frac{1}{\pi}\int_0^{\infty}\int_0^{\pi}\frac{1}{r}(a-a_r)\partial_{\theta}b d\theta dr-\frac{1}{\pi}\int_0^{\infty}\int_0^{\pi}\frac{1}{r}a_r\partial_{\theta}b d\theta dr  
\end{eqnarray}
We estimate the last two terms in \rec{bcterm}.
\par\medskip
{\bf Estimate of $\frac{1}{\pi}\int_0^{\infty}\int_0^{\pi}\frac{1}{r}(a-a_r)\partial_{\theta}b d\theta dr.$}
There holds:
\begin{eqnarray*}
\frac{1}{\pi}\int_0^{\infty}\int_0^{\pi}\frac{1}{r}(a-a_r)\partial_{\theta}b d\theta dr
&=&{ \frac{1}{\pi}\int\int_{R^2_+ }\nabla(\left\{\log((y_1-x_1)^2+y^2_2)^{1/2}) \right\})\cdot  [(a-a^+_{x_1}) \nabla^{\perp}b ]dy}\\
\end{eqnarray*}
Moreover we have the estimate:
\begin{eqnarray*}
\frac{1}{\pi}\int_0^{+\infty}\int_0^{\pi}\frac{1}{r}(a-a_r)\partial_{\theta} b dr\,d\theta
&\le&  \frac{1}{\pi}\left(\int_0^{+\infty} \int_0^\pi\frac{1}{r^2}|\partial_{\theta} a|^2 r dr\,d\theta\right)^{1/2}\left(\int_0^{+\infty} \int_0^\pi\frac{1}{r^2}|\partial_{\theta} b|^2 rdr\,d\theta\right)^{1/2}\\
&\leq& C\|\nabla a\|_{L^2(\R^2_+)}\|\nabla b\|_{L^2(\R^2_+).}
\end{eqnarray*}
\par\medskip
{\bf Estimate of $\frac{1}{\pi}\int_0^{\infty}\int_0^{\pi}\frac{1}{r}a_r\partial_{\theta}b d\theta dr.$}
\begin{eqnarray*}
\frac{1}{\pi}\int_0^{\infty}\int_0^{\pi}\frac{1}{r}a_r\partial_{\theta}b d\theta dr&=&\frac{1}{\pi}\int_0^{\infty}\frac{1}{r}a_r(b(r,\pi)-b(r,0)) dr\\&=&\frac{1}{\pi}\int_0^{\infty}\frac{1}{y_1}\left(\frac{a(y_1,0)+a(-y_1,0)}{2 }(b(-y_1,0)-b(y_1,0))\right)dy_1\\
&+&\frac{1}{2\pi}\int_{-\infty}^{\infty}\frac{1}{y_1}\left((a(y_1,0)+a(-y_1,0)) (b(-y_1,0)-b(y_1,0)\right)dy_1.
\end{eqnarray*}
Therefore the desired representation formula in $(0,0)$ is:
\begin{eqnarray}
w(0,0)&=&\frac{1}{\pi}\int\int_{\R^2_+ }\nabla( \log(|y|)  )\cdot  [(a-\frac{a(y_1,0)+a(-y_1,0)}{2}) \nabla^{\perp}b ]dy\nonumber\\ &
+&\label{eq:repform-0} \frac{1}{\pi}\int_{-\infty}^{\infty}\frac{1}{y_1}\left(\frac{(a(y_1,0)+a(-y_1,0))}{2} (b(-y_1,0)-b(y_1,0)\right)dy_1 \end{eqnarray}
\par
\medskip
For every $x_1\in\R$, $a^+_{x_1}(y_1):=\frac{a(x_1+y_1,0)+a(x_1 -y_1,0)}{2}$. 
We the get the representation formula for a generic point $(x_1,0)\in\partial\R^2_+$:
\begin{eqnarray}
w(x_1,0) 
&=&{ \frac{1}{\pi}\int\int_{\R^2_+ }\nabla(\left\{\log((y_1-x_1)^2+y^2_2)^{1/2}) \right\})\cdot  [(a-a^+_{x_1}) \nabla^{\perp}b ]dy}\nonumber\\
&+&{ \frac{1}{\pi}\int_{-\infty}^{\infty}\frac{1}{y_1} \left[a^+_{x_1}(y_1) (b(x_1-y_1,0)-b(x_1+y_1,0))\right]dy_1}
\end{eqnarray}
\par
{\bf Step 2.}  If  $a,b\in (H^{1} \cap L^{\infty})(\R^2_+)$ then we get the previous formula by approximation arguments.

\subsection{A Counter-Example to    $L^{\infty}$-Estimates }\label{CEx}
In this Section we will provide a counter-example to  Wente type estimates for the problem \rec{NPCompIntr}. Precisely we will show that even in the case  $a,b\in (H^1\cap L^{\infty})(\R^2_+)$ the solution given by
\eqref{reprform} needs not to be bounded.\par

Let $\psi\colon\R^2\to[0,+ \infty)$ be a radial smooth function such that:
\begin{equation}\label{psi}
\psi(x,y)=\left\{\begin{array}{cc}
1 &  (x,y)\in B(0,1/4),\\
0 & (x,y)\in B^c(0,1/2),\end{array}\right.\end{equation}
Let $\chi\colon\R\to \R$ be smooth  such that:
$$
\chi(x)=\left\{\begin{array}{cc}
1 &\text{if } x\ge 1,\\
0 &\text{if } x\le -1 \end{array}\right.$$
and $|\chi^{\prime}|\le C.$ Take for instance 
$$
\chi(x)=\left\{\begin{array}{cc}
\frac{2}{\pi}[\arctan(x)+\frac{\pi}{4}]& \text{if } -1\le x\le 1,\\
1 & \text{if }x\ge 1,\\
0 &\text{if } x\le -1, \end{array}\right.$$
We observe that $\chi(\frac{x}{\varepsilon})$ converges as $\varepsilon\to 0$ to the Heaviside function:
$$
H(x)=\left\{\begin{array}{cc}
1 &\text{if } x\ge 0,\\
0 & \text{if } x <0.\end{array}\right.$$
\begin{Prop}\label{example}
Let $\beta\in\R$ and consider the function:
\begin{equation}\label{b}
f(x)=(-\log|x|)^{-{\beta}}\psi(x).
\end{equation} 
Then \par
\begin{enumerate}
\item[$i)$] if $\beta\ge 0$, $f(x)\in ( H^{1/2}\cap L^{\infty}) (\R)$;\par
\item[$ii)$]   if $1/2<\beta$,  $f(x)H(x)\in (H^{1/2}\cap L^{\infty})(\R)$.\par
\end{enumerate}
\end{Prop}
{\bf Proof of Proposition \ref{example}.}
We prove only $ii)$. The proof of $i)$ is similar and even simpler.\par

It is clear that $f(x)H(x)\in L^{\infty}(\R).$ \par
$f(x)H(x)$ can be seen as the trace of the following function:
\begin{eqnarray*}  \tilde f(x,y)&=& (-1/2\log(x^2+y^2))^{-{\beta}}\psi(\sqrt{x^2+y^2}))\chi\left(\frac{x}{y}\right)  
\end{eqnarray*}

{\emph{Claim}:} $ \tilde f(x,y)\in H^{1}(\R^2_+),$ (this implies that $f(x)H(x)\in H^{1/2}(\R)$). \par
{\emph{Proof of the Claim}.} We estimate the $L^2$ norm of its partial derivatives.
\par\bigskip
{\bf Derivatives of $\tilde f:$}
\begin{eqnarray*}
\tilde f_y(x,y)&=&\chi^{\prime}\left(\frac{x}{y}\right)\left(-\frac{x}{y^2}\right)(-1/2\log(x^2+y^2))^{-{\beta}} \psi(\sqrt{x^2+y^2}))\\[5mm]
&+& \chi\left(\frac{x}{y}\right)\partial_y\left( (-1/2\log(x^2+y^2))^{-{\beta}}  \psi(\sqrt{x^2+y^2})\right)\\[5mm]
\tilde f_x(x,y)&=&\chi^{\prime}\left(\frac{x}{y}\right)\frac{1}{y}(-1/2\log(x^2+y^2))^{-{\beta}} \psi(\sqrt{x^2+y^2}))\\[5mm]
&+& \chi\left(\frac{x}{y}\right)\partial_x\left( (-1/2\log(x^2+y^2))^{-{\beta}}  \psi(\sqrt{x^2+y^2})\right).
\end{eqnarray*}

\medskip

{\bf $L^2$-estimate of $ \tilde f_y(x,y):$}
\par
\medskip
\begin{eqnarray*}
&&\int\int_{\R^2_+}|\tilde b_y(x,y)|^2dx dy\\
&&\lessapprox
\underbrace{ \int_{-1/2}^{1/2}\int_{|x|}^{1/2}(\chi^{\prime}\left(\frac{x}{y}\right)^2\left(\frac{x^2}{y^4}\right)(-1/2\log(x^2+y^2))^{-{(2\beta)}}   \psi^2(\sqrt{x^2+y^2}))\, dx dy}_{(1)}\\
&&+\underbrace{ \int\int_{(x^2+y^2)^{1/2}<1/2 } |\partial_y\left( (-1/2\log(x^2+y^2))^{-{\beta}} \psi(\sqrt{x^2+y^2})\right)|^2\, dx dy}_{(2)}.
\end{eqnarray*}
$\bullet$ We estimate $(2).$
\begin{eqnarray*} 
(2)&=&\int\int_{(x^2+y^2)^{1/2}<1/2 }\left|\psi^{\prime}(\sqrt{x^2+y^2})\frac{y}{(x^2+y^2)^{1/2}}(-1/2\log(x^2+y^2))^{-{\beta}}\right|.\\&&~~~~+\left.\psi(\sqrt{x^2+y^2}) \beta (-1/2\log(x^2+y^2))^{-{(\beta+1)}}\frac{2y}{x^2+y^2}\right|^2 dxdy\\
&\leq& C \int\int_{(x^2+y^2)^{1/2}<1/2 }(-1/2\log(x^2+y^2))^{-{2\beta}}\\&&~~~~~+(-1/2\log(x^2+y^2))^{-{2(\beta+1)}} \frac{4y^2}{(x^2+y^2)^2}dx dy<+\infty.
\end{eqnarray*}
\par

\noindent $\bullet$  We   estimate $(1)$, by recalling that  $\chi^{\prime}\left(\displaystyle\frac{x}{y}\right)\ne 0$ iff  $\displaystyle\frac{|x|}{|y|}\le 1$ and that

\begin{eqnarray*}&&\partial_y(-\frac{1}{y}(-1/2\log(x^2+y^2))^{-{\beta}})\\ 
&&=\frac{1}{y^2}(-1/2\log(x^2+y^2))^{-{\beta}})-\frac{1}{y}\partial_{y}((-1/2\log(x^2+y^2))^{-{\beta}})\\
&=&\frac{1}{y^2}(-1/2\log(x^2+y^2))^{-{\beta}})-\beta\frac{1}{y}\frac{y}{x^2+y^2}(-1/2\log(x^2+y^2))^{-{(1+\beta)}}).\end{eqnarray*}
We observe that
$$\beta\frac{1}{y}\frac{y}{x^2+y^2}(-1/2\log(x^2+y^2))^{-{(1+\beta)}})=o(\frac{1}{y^2}(-1/2\log(x^2+y^2))^{-{\beta}}))~~\mbox{as $(x,y)\to (0,0).$}$$
Therefore 
if $(x,y)\in B(0,1/2)$ we have 
$$
\frac{1}{y^2}(-1/2\log(x^2+y^2))^{-{\beta}})\le C \partial_y(-\frac{1}{y}(-1/2\log(x^2+y^2))^{-{\beta}}).$$
Hence:
\begin{eqnarray*}
(1)&=&
\int_{-1/2}^{1/2}\int_{|x|}^{1/2}(\chi^{\prime}\left(\frac{x}{y}\right)^2\left(\frac{x^2}{y^4}\right)(-1/2\log(x^2+y^2))^{-{2\beta}}   \psi^2(\sqrt{x^2+y^2}))\\
&\lessapprox
&\int_{-1/2}^{1/2}\int_{|x|}^{1/2}(\chi^{\prime}\left(\frac{x}{y}\right)^2\left(\frac{x^2}{y^4}\right)\left((-1/2\log(x^2+y^2))^{-{\beta}}\right) dx dy\\
&\le& \int_{-1/2}^{1/2}\int_{|x|}^{1/2}\partial_y(-\frac{1}{y}(-1/2\log(x^2+y^2))^{-{\beta}})\, dy dx\\
& \le& \int_{-1/2}^{1/2}\left[\frac{1}{|x|}(-1/2\log(x^2+x^2))^{-{\beta}}\right]-2(-1/2\log(x^2+\frac{1}{4}))^{-{\beta}} dx<+\infty.
\end{eqnarray*}

Observe that  since $\beta>1/2,$ the last integral is convergent.\par\medskip

{\bf $L^2$-estimate of $ \tilde f_x(x,y)$.}

\begin{eqnarray*}
&&\int\int_{\R^2_+}|\tilde f_x(x,y)|^2dx dy\\
&&\lessapprox
\underbrace{ \int_{-1/2}^{1/2}\int_{|x|}^{1/2}(\chi^{\prime}\left(\frac{x}{y}\right)^2\left(\frac{1}{y^2}\right)(-1/2\log(x^2+y^2))^{-{\beta}}   \psi^2(\sqrt{x^2+y^2}))\, dx dy}_{(3)}\\
&&+\underbrace{ \int\int_{(x^2+y^2)^{1/2}<1/2 } |\partial_x\left( (-1/2\log(x^2+y^2))^{-{\beta}} \psi(\sqrt{x^2+y^2})\right)|^2\, dx dy}_{(4)}.
\end{eqnarray*}
The estimate of $(3)$ is similar to $(1)$ and the estimate of $(4)$ is similar to the estimate of $(2)$.
\par
We can conclude the proof of Proposition \ref{example}.~~\hfill $\Box$

\subsection*{Estimate of $w(0,0).$}\par
Let us come back to the situation of subsection \ref{subsec:repform} and
consider:
$$a(x)=\psi(x)~~~\mbox{and}~~~b(x)=(-\log|x|)^{-{\beta}}\psi(x) H(x).$$
where $1/2<\beta<1$ and  $\psi$ is defined in \rec{psi}.\par
Since $b\equiv 0 $ in $y_1\le 0$ and $a$ is symmetric we have,
from \eqref{eq:repform-0}
\begin{eqnarray}
w(0,0)&=&\frac{1}{\pi}\int\int_{R^2_+ }\nabla(\left\{\log(|y|) \right\})\cdot  \left[\left(a-\frac{a(y_1,0)+a(-y_1,0)}{2}\right) \nabla^{\perp}b \right]dy\nonumber \\ 
&-&\frac{1}{\pi}\int_{0}^{+\infty}\frac{1}{y_1}\left(a(y_1,0) b(y_1,0)\right)dy_1 \end{eqnarray}
We already know that the first integral is finite. As for the second one,
 since we have chosen  $ 1/2<\beta<1$ we see that:
\begin{eqnarray*}
&&\frac{1}{\pi}\int_{0}^{+\infty}\frac{1}{y_1}\left((a(y_1,0) (b(y_1,0))\right)dy_1 =\frac{1}{\pi}\int_{0}^{+\infty}\frac{1}{y_1} (-\log|y_1|)^{-\beta}\psi^2(y_1)dy_1
\end{eqnarray*}
is divergent.
Hence $w(0,0)$ is not bounded.
\subsection{A Counter-Example to   $H^1$-estimates }

Consider now the solution of the   problem with vanishing Neumann boundary conditions:
\begin{equation}\left\{\begin{array}{cc}
-\Delta v_1=\nabla a\cdot\nabla^{\perp}b~~~\text{in }\R^2_+\\
\partial_{y_2}v_1=0~~~\text{on }\partial\R^2_+
\end{array}\right.
\end{equation}
given by the the representation formula:
$v_1(x)=\int_{D}\mathcal{G}(x,y)\nabla a(y)\cdot\nabla^{\perp}b(y)\, dy$.
By the same computations in subsection \ref{CEx} we find that: 
\begin{eqnarray}\label{v1}
v_1(x_1,0)&=&\frac{1}{\pi}\int\int_{R^2_+ }\nabla(\left\{\log((y_1-x_1)^2+x^2_2)^{1/2}) \right\})\cdot  [(a-a^+_{x_1}) \nabla^{\perp}b ]dy \\&+&
\frac{1}{2\pi}\int_{-\infty}^{\infty}\frac{1}{y_1} \left[\frac{(a(x_1+y_1,0)+a(x_1+y_1,0))}{2} (b(x_1-y_1,0)-b(x_1+y_1,0))\right]dy_1\nonumber\\
&-&{ \frac{1}{\pi}\int_{-\infty}^{+\infty}\log(|y_1-x_1|)[a\partial_{y_1} b ]dy_1}\nonumber
\end{eqnarray}
We take again:
\begin{equation} 
a(x)= \psi(x)~~\mbox{and}~~
b(x)=(-\log|x|)^{-{\beta}}\psi(x),
\end{equation}
with $0<\beta<1/2$ and $\psi$ defined as in \rec{psi}.
In this case  the solution $v_1$ is not in 
  $H^{1/2}(\R)$.
Indeed   if $0<\beta<1/2$, we have that:
$$a\partial_{x_1} b =\psi(x)[(-\beta)(-\log|x|)^{-{(\beta+1)}}\frac{1}{x}\psi(x)+\psi^{\prime}(x)(-\log|x|)^{-{\beta}}\notin  H^{-1/2}(\R)$$
One can check it by putting it in duality with $f(x)=[(\log|x|)^{-}]^{\beta}\in H^{1/2}(\R)$.

Now we observe that
in the representation of $v(x_1,0)$ the sum of the first two terms gives a function in ${H}^{1/2}(\R)$ ( it is the trace of a solution of the problem  Neumann problem \rec{pbhs} which is 
in $H^{1}(\R^2_+)$), the third term  

$${ \frac{1}{\pi}\int_{-\infty}^{+\infty}\log(|y_1-x_1|)[a\partial_{y_1} b ]dy_1}$$
cannot be in ${H}^{1/2}(\R)$ since $[a\partial_{x_1} b ] \notin  H^{-1/2}(\R).$ Therefore 
$v_1(x_1,0)\notin H^{1/2}(\R).$
\par
\bigskip
{\bf\small Acknowledgements.}
{\small The authors would like to thank Tristan Rivi\`ere  for  interesting discussions on the subject
and in particular for having pointed us out a mistake on the topic in the literature.
  While completing this work we heard of the existence of a work in preparation by Jonas Hirsch on similar questions.}

\end{document}